\documentclass[10pt english]{article}
\usepackage[english]{babel}
\usepackage[T1]{fontenc}
\usepackage{latexsym}
\usepackage{amssymb}
\usepackage{stmaryrd}
\usepackage{amsmath}
\usepackage{paralist} %liste numerotee
\usepackage[dvips]{graphicx,epsfig}
\usepackage[thmmarks, amsmath]{ntheorem}
\usepackage{color}

\theoremheaderfont{\sc}

\theoremheaderfont{\bf}
\theoremseparator{.}
\newtheorem{defi}{Definition}

\theoremseparator{.}
\newtheorem{thm}{Theorem}

\theoremseparator{.}

\theoremseparator{.}

\theoremseparator{.}
\newtheorem{lem}{Lemma}

\theoremseparator{.}
\theorembodyfont{\upshape}
\newtheorem{rem}{Remark}

\theoremseparator{.}
\theoremheaderfont{\it}
\theorembodyfont{\upshape}
\newtheorem{anex}{Example}

\theoremseparator{.}
\theoremheaderfont{\sl}
\theorembodyfont{\upshape}
\theoremprework{\relax}
\theorempostwork{\relax}
\theoremsymbol{\ensuremath{\Box}}
\newtheorem*{pro}{Proof}

\theoremseparator{.}
\theoremheaderfont{\sc}
\theorembodyfont{\upshape}
\theoremprework{\relax}
\theorempostwork{\relax}
\theoremsymbol{\ensuremath{\Box}}

\theoremseparator{}
\theoremheaderfont{\sc}
\theorembodyfont{\upshape}
\theoremprework{\relax}
\theorempostwork{\relax}
\theoremsymbol{\ensuremath{\Box}}

\newcommand{\shs}{\hspace{0.2cm}}
\newcommand{\hs}{\hspace{0.5cm}}

\begin{document}

\noindent\textbf{Title:} Ergodic properties of  geometrical crystallization processes, I\\
\noindent\textbf{Authors:} Y. Davydov (University of Lille 1) and A. Illig (University of Versailles Saint-Quentin)\\
\noindent\textbf{Subj-class:} PR\\
\noindent\textbf{MSC-class:} 60G60 (primary), 37A25 (secondary)\\
\noindent\textbf{Abstract:} We are interested here in a birth-and-growth process where germs are born according to a Poisson point process with invariant under translation in space intensity measure. The germs can be born in free space and then start growing until occupying the available space. In order to consider various way of growing, we describe the crystals at each time through their geometrical properties. In this general framework, the crystallization process can be caracterized by the random field giving for a point in the space state the first time this point is reached by a crystal. We prove under general conditions that this random field is mixing in the sens of ergodic theory and obtain estimates for the coefficient of absolute regularity.\\


%------------------------------------------------------------------------------
\section{Introduction}
%------------------------------------------------------------------------------
The crystallization process we consider here deals with germs $g=(x_g,t_g)$ that appear at random times $t_g$ on random locations $x_g$. The born process $\mathcal{N}$ is a Poisson point process on $\mathbb{R}^d\times\mathbb{R}^+$ with intensity measure denoted by $\Lambda$. Once the germs or cristallization centers are born, crystals are allowed to grow if their location is not still occupied by another crystal and when two crystals meet the growth stops a meeting points. There are then many ways to describe crystal expansion. The first approach is to consider the random sets (called crystallization state) that corresponds to the fraction of space occupied by crystals at a given time. In this case, crystallization is studied through the theory of set-valued processes. Another way to describe crystal growth is to deduce  the expression of the speed growth  from characteristic local media properties of the state space. One can also consider for a germ $g\in\mathbb{R}^d\times\mathbb{R}^+$ and a point $x\in\mathbb{R}^d$ the time $A_g(x)$ at which $x$ is reached by the free crystal associated to the germ $g$. The crystallization process is then caracterized by the following random field $\xi$ giving for a location $x\in\mathbb{R}^d$ its crystallization time $\xi(x)$:
\begin{equation}\label{process}\xi(x)=\inf_{g\in\mathcal{N}}A_g(x).
\end{equation}
We adopt in this paper the last definition and study the crystallization process througth the random field $\xi$. 

This model was introduced by Kolmogorov \cite{Kol37} and  independently
by Johnson \& Mehl \cite{JM39}, and intensively studied by many authors. We mention here only a few number of papers which represent the main approaches
and where one can find an exhaustive liste of references :
M\o ller \cite{Mol89}, \cite{Mol92}, and also Micheletti \& Capasso \cite{MC97}.
A very large part of these investigations deals with the geometrical structure of the mosaic once all the germs have finished their growth. Here, we  are rather interested in estimation problems (such as the estimation of the parameters of the intensity measure $\Lambda$ or other functionals like the number of crystals in the limit mosaic) in the case when only one realisation can be observed on a sufficiently large domain compared to the mean size of crystals. Naturally, we suppose that the crystallisation process is space homogeneous. More precisely, we assume that the intensity measure is defined as follows,
\begin{equation}\label{product}
\Lambda=\lambda^d\times m,
\end{equation}
where $\lambda^d$ is the Lebesgue measure on $\mathbb{R}^d$ and $m$ is a measure on $\mathbb{R}^+$ finite on bounded Borelians.

This article is mainly devoted to ergodic properties of the random field $\xi(x)_{x\in \mathbb{R}^d}$ defined by (\ref{process}) which deliver a solid base for efficient estimation of parameters of the model and subsequent application to the study of its asymptotical normality. Under the above hypothesis and rather general conditions on growth speed and geometrical shape of crystals, we demonstrate that the random field $\xi$ is mixing in the sens of the ergodic theory. Moreover, under some additional assumptions, we obtain estimates of the absolute regularity coefficient of $\xi$. 

The statistical application represents de second part of our work and will be published elsewhere.
%The last section deals with estimation problem on examples.

%------------------------------------------------------------------------------
\section{Assumptions on the birth-and-growth process}
%------------------------------------------------------------------------------
\subsection{The birth process}
Germs are born according to a Poisson point process on $E=\mathbb{R}^d\times \mathbb{R}^+$ denoted by $\mathcal{N}$. Thus, a germ is a random point $g=(x_g,t_g)$ in $\mathbb{R}^d\times \mathbb{R}^+$, where $x_g$ is the location in the growing space  $\mathbb{R}^d$ and $t_g$ is the time of birth on the time axes $\mathbb{R}^+$. We suppose that the intensity measure $\Lambda$ of $\mathcal{N}$ is the product (\ref{product}) of the Lebesgue measure $\lambda^d$ on $\mathbb{R}^d$ and a measure $m$ on $\mathbb{R}^+$ such that $m([0,a])<\infty$ for all $a>0$.  The most interesting cases to be considered (see \cite{Mol86}) are for a discret measure $m$ or when $m(dt)=\alpha t^{\beta-1}\lambda(dt)$ with $\alpha>0$ and $\beta>0$.
 
Since the Lebesgue measure is invariant by the translations on $\mathbb{R}^d$, 
we derive that $\mathcal{N}$ is space homogeneous. So we are led to consider only sets around the origine. In particular, we introduce for a time $t$, the causal cone:
\begin{equation}\label{cone}
K_t=\{g\in E\,|\,A_g(0)\leq t\}
\end{equation}
which consists of all the possible germs that are able to reach the origine before time~$t$. The measure $\Lambda(K_t)$ of the causal cone $K_t$ is denoted by $F(t)$. 
\subsection{Expansion of crystals}
We call ``free crystal'' a crystal which is born in a fraction of space non-occupied by other crystals at the time of its birth. We associate to each germ $g$ in $E$ a function $A_g$:
\begin{equation}
\label{fonction}
\begin{array}{llcl}
A_g: & \mathbb{R}^d & \rightarrow & \mathbb{R}^+\\
 & x & \mapsto & A_g(x)
\end{array}
\end{equation}
where $A_g(x)$ is the time when $x$ is reached by the crystal assumed to be free and associated to the germ $g$. Consequently, at time $t$ a free crystal is defined by the set 
\begin{equation}\label{crystal}
C_g(t)=\{x\,|\,A_g(x)\leq t\}.
\end{equation}

In the following, we make several assumptions on the free crystals family $\{C_g,\;g\in\mathcal{N}\}$ and the functions family $\{A_g,\;g\in\mathcal{N}\}$. We also specify when necessary the link between assumptions and crystal growth.

\begin{enumerate}[1)]
\item $\forall\, g=(x_g,t_g)\in E,\shs A_g(x)=A_{(0,t_g)}(x-x_g)\shs\forall\,x\in\mathbb{R}^d.$

Crystal growth is space homogeneous. This assumption implies that for all germ $g=(x_g,t_g)$,
$$C_g(t)=C_{(0,t_g)}(t)+x_g\hs\forall \,t\in\mathbb{R}^+.$$

\item $\forall\, g=(x_g,t_g)\in E,\shs A_g(x_g)=t_g\textrm{ and }A_g(x)\geq t_g\shs\forall \,x\in\mathbb{R}^d.$

A crystal can only reach a point $x$ after its birth.

\item The free crystals $C_g(t)$ are bounded, convex sets and the family $(C_g(t))_{t\in \mathbb{R}^+}$ is increasing that means
$$C_g(s)\subset C_g(t)\hs\forall \,0\leq s<t.$$

\item The functions $x\mapsto A_g(x)$ are continuous.

Thus, crystals grow in each space direction and without any jump so that
$$\partial C_g(t)=\{x\,|\,A_g(x)=t\}.$$

\item There exists $M>0$ such that $\forall\,t_g\in\mathbb{R}^+ $ we have
$$A_{(0,t_g)}(x)\geq t_g+\frac{1}{M}|x|\hs\forall \,x\in\mathbb{R}^d.$$
The growth speed is then bounded by the constant $M$.

\item $\forall \,g=(0,t_g), \shs\forall \,r>0, \shs \exists \,t>0$ such that 
$$C_g(t)\supset B(0,r)=\{x\in\mathbb{R}^d\,|\,|x|\leq r\}.$$
A free crystal grows in each direction and never definitively stops growing.

\item If $L_g=\{(x,t)\,|\,A_g(x)\leq t\}$ denotes the epigraph of $A_g$, then $\forall \,g_1\in L_g$,
$$A_g(x)\leq A_{g_1}(x)\hs\forall \,x\in\mathbb{R}^d.$$
This means that a crystal born inside $L_g$ never exits.
\end{enumerate}
When $d=1$, we introduce for each germ $g$ the restrictions  $A^+_g$, $A^-_g$ of $A_g$ respectively to $[x_g,+\infty)$, $(-\infty,x_g]$ and consider when necessary a stronger version of Assumption $7)$:
\begin{enumerate}[7a)]
\item $\forall g_1 \in E,\,\forall g_2\in E$, if for some $x_0\geq x_g$, $A^+_{g_1}(x_0)\geq A_{g_2}(x_0)$ then 
$$A_{g_1}(x)\geq A_{g_2}(x) \hs\forall x\geq x_0$$
and if for $x_1\leq x_g$, $A_{g_1}^-(x_1)\geq A_{g_2}(x_1)$ then
$$A_{g_1}(x)\geq A_{g_2}(x) \hs\forall x\leq x_1.$$
\end{enumerate}
Some remarks can be made on a part of these assumptions.

\begin{rem} The assumption $3)$ implies that for all $0<s<t$, 
$$C_g(s)\subset[C_g(t)]^{\circ}.$$
Indeed, if there exists $x\in \partial C_g(t)\cap C_g(s)$, then we should have$A_g(x)=t$ and $A_g(x)\leq s$. But, this cannot occur.
\end{rem}

\begin{rem} Let $g=(0,t_g)$ be a germ. Observe that,
$$\sup_{|x|\leq r}A_g(x)=\sup_{|x|=r}A_g(x).$$
If $x_g$ satisfies $|x_g|=r$ and $t_g=A_g(x_g)=\sup_{|x|= r}A_g(x)$ then $C_g(t_g)$ contain all the points $x$ such that $|x|=r$ and by convexity (assumption $3)$), we obtain that $C_g(t_g)\supset B(0,r)$.
\end{rem}

\begin{rem} \label{Mg}The function $r\mapsto M_g(r)=\sup_{|x|=r}A_g(x)$ is increasing on $\mathbb{R}^+$ and 
$$\lim_{r\rightarrow\infty}M_g(r)=+\infty.$$

\noindent To prove the first assertion, let us consider $0\leq r_0<r_1$ and $(x_g,t_g)\in\mathbb{R}^d\times\mathbb{R}^+$ such that $|x_g|=r_0$ and $t_g=A_g(x_g)=M_g(r_0)$. Then, $x_g\in\partial C_g(t_g)$ and by convexity, $C_g(t_g)\supset B(0,r_0)$. If $M_g(r_1)=M_g(r_0)$, $B(0,r_1)$ would also be included in $C_g(t_g)$ and $x_g$ would be inside $[C_g(t_g)]^{\circ}$. But, this is impossible.\\
For the second point, the assumption $5)$ implies
 that
$$A_{(0,t_g)}(x)\geq t_g+\frac{1}{M}|x|$$
and
$$M_g(r)\geq t_g+\frac{r}{M}.$$
Thus, $\lim_{r\rightarrow\infty}M_g(r)=+\infty$.
\end{rem}

To obtain the absolute regularity property of $\xi$ when $d$ is greater or equal to $2$, we add two other assumptions. Let us introduce some definitions before stating the assumptions. For a germ $g=(x_g,t_g)$ and the associated  free crystal $C_g(t)$ at time $t$, we call ``interior diameter'' and write $d_g(t)$ the diameter of the greatest ball centered in $x_g$ and included in $C_g(t)$. In the same way, $D_g(t)$ named the ``exterior diameter'' denotes the diameter of the smallest ball centered in $x_g$ containing $C_g(t)$. From the preceding assumptions, we deduce that for any germe $g\in E$, the functions $t\mapsto d_g(t)$ and $t\mapsto D_g(t)$ are continuous and for all $t\in\mathbb{R}^+$, we have $d_g(t)\leq D_g(t)$. The additional assumptions are the following ones:
  
\begin{enumerate}[8)]
\item $\exists A>0$ such that $\forall\, g\in E,\shs\forall\, t\in\mathbb{R}^+$,
$$\frac{1}{A}D_g(t)\leq d_g(t).$$
This assumption ensure that free crystals have non-degenerated shapes.
\end{enumerate}

\begin{enumerate}[9)]
\item $\forall g=(x_g,t_g)\in E$, the function $t\mapsto D_g(t)$ is ``subadditive'':
$$D_g(t+h)\leq D_g(t)+D_{(0,t)}(h)\hs\forall\, t\geq t_g,\shs\forall\, h\geq 0.$$
\end{enumerate}

We give now an example that satisfies all the assumptions from $1)$ to $9)$.

\begin{anex} \label{modelsimple}For any germ $g=(x_g,t_g)$, we suppose that the crystal at time $t\geq t_g$ is as follows:
$$C_g(t)=x_g+[V(t)-V(t_g)]K,$$ 
where $K$ is a convex compact set such that $0\in K^{\circ}$ and the function $t\mapsto V(t)$ represents the distance achieved with function speed $t\mapsto v(t)$. We assume that $v$ is positive almost everywhere. Moreover, we suppose that $V$ is absolutely continuous:
$$V(t)=\int_{0}^t v(s)ds\hs\forall t\in\mathbb{R}^+$$
 and such that for all $t\geq 0$, $h>0$,
$$V(t+h)-V(t)>0.$$
Observe that
$$C_g(t+h)=C_g(t)\varoplus [V(t+h)-V(t)]K$$
where $\varoplus$ represents here the Minkowski summation of two sets $A$ and $B$:
$$A\varoplus B=\{x+y\,|\,x\in A,\;y\in B\}.$$

Now, we denote by $p_{x,K}$ the norm of the intersection point between $\partial K$ and the line $(0,x)$. Then, a point $x$ is reached at time $t$ by the crystal born in $x_g$ at time $t_g$ if
\begin{equation}\label{ass5}
(V(t)-V(t_g))p_{x-x_g,K}=|x-x_g|.
\end{equation}
%%%%%%%%%%%%%%%%%%%%%
As $V$ is invertible, 
$$t=A_g(x)=V^{-1}\left(\frac{|x-x_g|}{p_{x-x_g,K}}+V(t_g)\right).$$
Thus, all the assumptions $1)$ to $9)$ except assumption $5)$ are satisfied in this example. 
%%%%%%%%%%%%%%%%%%%
For the last assumption, we can suppose for example that $v$ is bounded, $v(s) \leq L.$ 

We can take $x_g = 0$. As  $K$ is compact, there exists a constant $C$ such that
$p_{x,K} \leq C$ for all $x$. From (\ref{ass5}) we get
$$
|x| = (V(t)-V(t_g))p_{x,K}\leq LC_{}(A_{g}(x)-t_{g}),
$$
which gives assumption 5) with the constant $M = LC_{}.$

%%%%%%%%%%%%%%%%%%%%

Note that if $K=B(0,1)$ and $v(t)=c$, this example corresponds to the linear homogeneous expansion in all directions.
\end{anex}

%------------------------------------------------------------------------------
\section{Mixing property}
%------------------------------------------------------------------------------

We assume without loss of generality that the random field  $\xi=(\xi(x))_{x\in \mathbb{R}^d}$ defined by (\ref{process}) is a canonical random field on $(\Omega,\mathcal{F},\mathbb{P})$. Namely, we suppose that $\Omega=\mathbb{R}^{T}$ with $T=\mathbb{R}^d$, $\mathcal{F}$ is the $\sigma$-algebra generated by the cylinders and $\mathbb{P}$ is the distribution of $\xi$ so that for all $\omega\in\Omega$, $\xi(x,\omega)=\omega(x)$. As Lebesque measure $\lambda^d$ on $\mathbb{R}^d$ is invariant by the translations on $\mathbb{R}^d$, we deduce that $\xi$ is homogeneous. This means that $\mathbb{P}$ is invariant by the translations 
$$S_h(\omega)(x)=\omega(x+h),\;h\in\mathbb{R}^d.$$
We precise here what we call a mixing random field.
\begin{defi}\label{mixdef} A random field $\xi=(\xi(x))_{x\in \mathbb{R}^d}$ is mixing if for all $A$, $B\in\mathcal{F}$,
\begin{equation}\label{cyl}
\mathbb{P}(A\cap S^{-1}_h(B))\xrightarrow[|h|\rightarrow\infty]{}\mathbb{P}(A)\mathbb{P}(B).
\end{equation}
\end{defi}
\begin{rem} We note here that if a random field is mixing in the sense of Definition \ref{mixdef}, then the random field is also egodic.
\end{rem}
To prove that a random field is mixing, it is sufficient to verify Condition (\ref{cyl}) for cylinders and establish the following condition
\begin{equation}\label{cond} 
\begin{split}
 & \forall x_1,\dots,x_k\in\mathbb{R}^d,\;\forall y_1,\dots,y_m\in\mathbb{R}^d,\;\forall E_1\in\mathcal{B}^k,\;\forall E_2\in\mathcal{B}^m,\\
& \mathbb{P}\left\{(\xi(x_1),\dots,\xi(x_k))\in E_1,\;(\xi(y_1),\dots,\xi(y_m))\in E_2\right\}  \\
 & \xrightarrow[|h|\rightarrow\infty]{}  \mathbb{P}\left\{(\xi(x_1),\dots,\xi(x_k))\in E_1\right\}\mathbb{P}\left\{(\xi(y_1),\dots,\xi(y_m))\in E_2\right\}\\
\end{split}
\end{equation}
where $\mathcal{B}^k$ (respectively $\mathcal{B}^m$) is the $k$-dimentional (respectively $m$-dimentional) Borelian $\sigma$-field.
%---------------------------Theoreme 1-----------------------------------------
\begin{thm} \label{mixing} $(d\geq 1)$ Under assumtions $1)$ to $7)$, the random field  $\xi=(\xi(x))_{x\in \mathbb{R}^d}$ defined by (\ref{process}) is mixing.
\end{thm}
To demonstrate Theorem \ref{mixing}, we need three auxiliary lemmas.
%---------------------------- Lemme 1------------------------------------------
\begin{lem}\label{inegalite} If $A_1$, $A_2$, $B_1$ and $B_2$, are four events, then
\begin{itemize}
\item[(i)] $|\mathbb{P}(A_1)-\mathbb{P}(A_2)|\leq\mathbb{P}(A_1\triangle A_2)$,
\item[(ii)]$|\mathbb{P}(A_1\cap B_1)-\mathbb{P}(A_2\cap B_2)|\leq\mathbb{P}(A_1\triangle A_2)+\mathbb{P}(B_1\triangle B_2)$,
\end{itemize}
where for two events $A$ and $B$, $A\triangle B=(A\cap B^c)\cup(A^c\cap B)$.
\end{lem}
\begin{pro} Elementary.
\end{pro}
%----------------------------------Lemme 2 ------------------------------------
Now, for all $h\in\mathbb{R}^d$ and $r\geq 0$, we define new random fields to approximate $\xi(\cdot)$ and its translations $\xi(\cdot+h)$:
$$\xi_r^h(x)=\inf_{
\scriptsize
\begin{array}{c}
g\in\mathcal{N}\\
|x_g-h|\leq r
\end{array}}A_g(x)\hs\forall\, x\in\mathbb{R}^d.$$
\begin{lem} \label{lemme}Let $H(R)=M_{(0,R)}(R)$ with $M_g(r)$ defined in Remark \ref{Mg}. 
Under Assumptions $1)$ to $7)$, we have for all $h\in\mathbb{R}^d$,
$$\mathbb{P}\left(\xi(x+h)=\xi_{(M+1)H(R)}^h(x),\;|x|\leq R\right)\geq 1-\textrm{e}^{-F(R)}$$
where $M$ is the constant of Assumption $5)$ and $F(R)$ is the measure of the causal cone $K_R$ defined by relation (\ref{cone}).
\end{lem}
\begin{pro} As $\mathcal{N}$ is space homogeneous, 
$$\mathbb{P}\left(\xi(x+h)=\xi_{(M+1)H(R)}^h(x),\;|x|\leq R\right)=\mathbb{P}\left(\xi(x)=\xi_{(M+1)H(R)}^0(x),\;|x|\leq R\right)$$
and it is then sufficient to demonstrate Lemma \ref{lemme} for $h=0$. First, observe that Assumption $7)$ implies that
$$\{\xi(0)\leq R\}\subset\left\{\sup_{|x|\leq R}\xi(x)\leq H(R)\right\}.$$
Now, let us prove that
\begin{equation}\label{rel}
\left\{\sup_{|x|\leq R}\xi(x)\leq H(R)\right\}\subset\{\xi(x)=\xi_{R+MH(R)}^0(x),\;|x|\leq R\}.
\end{equation}
To prove (\ref{rel}), note that Assumptions $1)$ and $5)$ imply that for all germ $g$,
$$A_g(x)\geq t_g+\frac{|x-x_g|}{M}\hs\forall\, x\in\mathbb{R}^d.$$
In particular, for germs $g$ such that $|x_g|>R+MH(R)$, we deduce that
$$A_g(x)>H(R)\hs\forall\, x\in\mathbb{R}^d,\shs|x|\leq R.$$
Hence, for all $x$ such that $|x|\leq R$,
$$\inf_{\scriptsize
\begin{array}{c}
g\in\mathcal{N}\\
|x_g|> R+MH(R)
\end{array}}A_g(x)>H(R)\geq \xi(x)$$
and (\ref{rel}) follows. On the other and, for $0\leq r_1\leq r_2$
$$\xi(x)\leq \xi_{r_2}^0(x)\leq\xi_{r_1}^0(x)\hs\forall\, x\in\mathbb{R}^d.$$
From Assumption $2)$, remark that $R\leq H(R)$ and deduce that
$$\{\xi(x)=\xi_{R+MH(R)}(x),\;|x|\leq R\}\subset\{\xi(x)=\xi_{(M+1)H(R)}(x),\;|x|\leq R\}.$$
Finally, we obtain that
$$\mathbb{P}\left\{\xi(x)=\xi_{(M+1)H(R)}(x),\;|x|\leq R\right\}\geq\mathbb{P}\left\{\xi(0)\leq R\right\}$$
and 
$$\mathbb{P}\left\{\xi(0)\leq R\right\} = \mathbb{P}\left\{\mathcal{N}\cap K_{R}\ne\emptyset\right\}.$$
But, 
\begin{equation*}
\mathbb{P}\left\{\mathcal{N}\cap K_{R}\ne\emptyset\right\}=1-\textrm{e}^{-\Lambda(K_{R})}.
\end{equation*}
\end{pro}
%----------------------------Lemme 3-------------------------------------------
\begin{lem}\label{infinity}Assumptions $1)$ and $6)$ imply that
$$F(R)=\Lambda(K_{R})\xrightarrow[R\rightarrow\infty]{}\infty.$$
\end{lem}
\begin{pro} The assumptions $1)$ and $6)$ imply that for all germ $g\in E$, there exists $R>0$ such that $g\in K_{R}$ or equivalently such that $0$ belongs to the crystal $C_g(R)$. But, 
$$\bigcup_{R\geq 0} K_{R}=E$$
and since $\Lambda(E)=+\infty$, the result follows.
\end{pro}
We come back to the demonstration of Theorem \ref{mixing}.
\begin{pro}  For $(x_1,\dots,x_k)$ in $E^k$, $(y_1,\dots,y_m)$ in $E^m$, $E_1$ in $\mathcal{B}^k$ and $E_2\in\mathcal{B}^m$, we define the sets:
$$\begin{array}{lll}
A & = &\{(\xi(x_1),\dots,\xi(x_k))\in E_1\},\\
 & & \\
B & = &\{(\xi(y_1),\dots,\xi(y_m))\in E_2\},\\
 & & \\
B_h & = & \{(\xi(y_1+h),\dots,\xi(y_m+h))\in E_2\}.\\
\end{array}$$
Let us define $r=\max\{|x_i|,\;i=1\dots k;\;|y_j|,\;j=1\dots m\}$ and consider a positive real number  $\epsilon$. From Lemma \ref{lemme} and Lemma \ref{infinity}, we find $R>r$ such that 
$$\mathbb{P}\left\{\xi(x)=\xi^0_{(M+1)H(R)}(x),\;|x|\leq R\right\}\geq 1-\epsilon.$$
Let us now introduce $h\in\mathbb{R}^d$ such that $|h|>2R_1$, where $R_1=(M+1)H(R)$. We also define some other sets:
$$\begin{array}{lll}
\tilde A & = &\{(\xi_{R_1}^0(x_1),\dots,\xi_{R_1}^0(x_k))\in E_1\},\\
 & & \\
\tilde B & = &\{(\xi_{R_1}^0(y_1),\dots,\xi_{R_1}^0(y_m))\in E_2\},\\
 & & \\
\tilde B_h & = & \{(\xi_{R_1}^h(y_1),\dots,\xi_{R_1}^h(y_m))\in E_2\}.\\
\end{array}$$
Lemma \ref{inegalite} $(ii)$ leads to the following inequality:
$$|\mathbb{P}(A\cap B_h)-\mathbb{P}(\tilde A\cap\tilde B_h)|\leq\mathbb{P}(A\triangle\tilde A)+\mathbb{P}(B_h\triangle\tilde B_h).$$
We introduce the set $D=\{\xi(x)=\xi_{R_1}^0(x),\;|x|\leq R\}$ and obtain by Lemma \ref{lemme} that
$$\begin{array}{lll}
\mathbb{P}(A\triangle\tilde A) & = & \mathbb{P}((A\triangle\tilde A)\cap D)+\mathbb{P}((A\triangle\tilde A)\cap D^c)\\
 & = & \mathbb{P}((A\triangle\tilde A)\cap D^c)\\
 & \leq & \mathbb{P}(D^c)\\
 & \leq & \epsilon.
\end{array}$$
If we introduce the set $D_h=\{\xi(x+h)=\xi_{R_1}^h(x),\;|x|\leq R\}$ in place of $D$, we obtain by the same arguments that
$$\mathbb{P}(B_h\triangle\tilde B_h)\leq \epsilon.$$
These two inequalities imply that
\begin{equation}\label{in1}
|\mathbb{P}(A\cap B_h)-\mathbb{P}(\tilde A\cap\tilde B_h)|\leq 2\epsilon.
\end{equation}
On the other hand, the events $\tilde A$ and $\tilde B_h$ are independent because $|h|>2R_1$. Thus,
$$\mathbb{P}(\tilde A\cap \tilde B_h)=\mathbb{P}(\tilde A)\mathbb{P}(\tilde B_h)$$
and by space homogeneity of $\mathcal{N}$, $\mathbb{P}(\tilde B_h)=\mathbb{P}(\tilde B)$ so that
\begin{equation}\label{in2}
\mathbb{P}(\tilde A\cap \tilde B_h)=\mathbb{P}(\tilde A)\mathbb{P}(\tilde B).
\end{equation}
Moreover, by Lemma \ref{inegalite} $(i)$ 
$$\begin{array}{lll}
|\mathbb{P}(A)\mathbb{P}(B)-\mathbb{P}(\tilde A)\mathbb{P}(\tilde B)| & \leq & |\mathbb{P}(A)-\mathbb{P}(\tilde A)|+|\mathbb{P}(B)-\mathbb{P}(\tilde B)|\\
 &  \leq & \mathbb{P}(A\triangle\tilde A)+\mathbb{P}(B\triangle\tilde B).
\end{array}$$
and Lemma \ref{lemme} implies that
\begin{equation}\label{in3}
|\mathbb{P}(A)\mathbb{P}(B)-\mathbb{P}(\tilde A)\mathbb{P}(\tilde B)|\leq 2\epsilon.
\end{equation}
Inequalities (\ref{in1}), (\ref{in3}) and relation (\ref{in2}) imply that for all $h\in\mathbb{R}^d$ such that $|h|>2R_1$,
$$
|\mathbb{P}(A\cap B_h)-\mathbb{P}(A)\mathbb{P}(B)|\leq 4\epsilon
$$ 
and Theorem \ref{mixing} is then proved.
\end{pro}

%------------------------------------------------------------------------------
\section{Absolute regularity}
%------------------------------------------------------------------------------

\subsection{General definitions}
%-------------------------------------------------------------
For a subset $T$ of $\mathbb{R}^d$, we denote by $\mathcal{F}_T$ the $\sigma$-field generated by the random variables $\xi(x)$ for all $x$ in $T$.  Now, consider two disjoint sets $T_1$ and $T_2$ in $\mathbb{R}^d$ and define the absolute regularity coefficient for the $\sigma$-fields $\mathcal{F}_{T_1}$ and $\mathcal{F}_{T_2}$ as follows:
$$\beta(T_1,T_2)=\|\mathcal{P}_{T_1\cup T_2}-\mathcal{P}_{T_1}\times\mathcal{P}_{T_2}\|_{var},$$
where $\|\mu\|_{var}$ is the total variation norm of a signed measure $\mu$ and $\mathcal{P}_T$ is the distribution of the restriction $\xi_{|T}$ in the set $\mathcal{C}(T)$ of continuous real-valued functions defined on T. If $T_1\cap T_2=\emptyset$, note that $\mathcal{C}(T_1\cup T_2)$ is canonically identified to $\mathcal{C}(T_1)\times\mathcal{C}(T_2)$. 

The strong mixing coefficient is defined as follows,
$$\alpha(T_1,T_2)=\sup_{
\begin{array}{c}
A\in\mathcal{F}_{T_{1}}\\
B\in\mathcal{F}_{T_{2}}
\end{array}}|\mathbb{P}(A\cap B)-\mathbb{P}(A)\mathbb{P}(B)|$$

The process $\xi$ is said to be absolutely regular ($\alpha$-mixing) if the absolute regularity coefficient (the strong mixing coefficient) converges to zero when the distance between $T_1$ and $T_2$ tends to infinity with $T_1$ and $T_2$ belonging to a certain class of sets.

\begin{rem}It is well known that
$$\alpha(T_1,T_2)\leq \frac{1}{2}\beta(T_1,T_2)$$
so that absolute regularity of the process $\xi$  implies $\alpha$-mixing.
\end{rem}

When $d=1$, one usually chooses $T_1=(-\infty,0]$ and $T_2=[r,+\infty)$ whereas in the case $d\geq 2$, there are several sorts of sets to be considered. The results we obtain in this paper when $d\geq 2$ deal with quadrant domains as represented on Figure $1$ and enclosed cube domains as represented on Figure $2$.

\subsection{Upper bounds}
%-----------------------
\subsubsection{Approach}
In order to obtain upper bounds for the absolute regularity coefficient $\beta(T_1,T_2)$, we approximate  the restrictions of $\xi$ on $T_1$ and $T_2$ by two independent random fields and apply the following lemma.
\begin{lem} \label{outil} Let us consider a random field $(\eta(x))_{x\in\mathbb{R}^d}$ and two disjoint subsets $T_1$ and $T_2$ of $\mathbb{R}^d$. If there exists two random fields $(\eta_1(x))_{x\in\mathbb{R}^d}$ and $(\eta_2(x))_{x\in\mathbb{R}^d}$ and two positive constants $\delta_1$ and $\delta_2$ such that:
\begin{itemize}
\item $\eta_1$ and $\eta_2$ are independent
\item $\mathbb{P}\left\{\xi(x)=\eta_i(x),\shs\forall\,x\in T_i\right\}\geq 1-\delta_i$ for $i=1,2$.
\end{itemize}
then
$$\beta(T_1,T_2)\leq 8\,(\delta_1+\delta_2).$$
\end{lem}

\begin{pro}Let us denote by $\mathcal{P}_1$ the distribution of the restriction $\xi_{|T_1}$ of $\xi$ to $T_1$, by  $\mathcal{P}_2$ the distribution of the restriction $\xi_{|T_2}$ of $\xi$ to $T_2$, by $\mathcal{Q}_1$ the distribution the restriction ${\eta_1}_{|T_1}$ of $\eta_1$ to $T_1$, and by $\mathcal{Q}_2$ the distribution of the restriction ${\eta_2}_{|T_2}$ of $\eta_2$ to $T_2$. We have for $i=1,2$, that 
$$\|\mathcal{P}_i-\mathcal{Q}_i\|_{var}\leq 4\delta_i.$$
Indeed, it is clear that
$$\|\mathcal{P}_i-\mathcal{Q}_i\|_{var}=2\,\sup_{A}|\mathcal{P}_i(A)-\mathcal{Q}_i(A)|.$$
If we denote by $D_i$ the set $\{\xi(x)=\eta_i(x),\;\forall\,x\in T_i\}$, we obtain that
$$\|\mathcal{P}_i-\mathcal{Q}_i\|_{var}=2\,\sup_A|\mathbb{P}(\{\xi_{|T_i}\in A\}\cap D_i^c)-\mathbb{P}(\{{\eta_i}_{|T_i}\in A\}\cap D_i^c)|$$
and deduce that
$$\|\mathcal{P}_i-\mathcal{Q}_i\|_{var}\leq 4\,\mathbb{P}(D_i^c).$$
Since $\mathbb{P}(D_i)\geq 1-\delta_i$, we conclude that
$$\|\mathcal{P}_i-\mathcal{Q}_i\|_{var}\leq 4 \,\delta_i.$$

Now, we denote by $\mathcal{P}$ the distribution of $\xi$ on $T_1\cup T_2$ and $\mathcal{Q}$ the disctribution of $\eta$ on $T_1\cup T_2$ with $\eta$ defined as follows:
$$\eta(x)=\left\{\begin{array}{ll}
\eta_1(x) & x\in T_1,\\
\eta_2(x) & x \in T_2.
\end{array}\right.$$
We have
$$\mathbb{P}\left\{\xi(x)=\eta(x),\shs \forall\, x\in T_1\cup T_2\right\}=\mathbb{P}(D_1\cap D_2).$$
But,
$$\mathbb{P}(D_1\cap D_2)=1-\mathcal{P}(D_1^c\cup D_2^c)\geq 1-\mathbb{P}(D_1^c)-\mathbb{P}(D_2^c)$$
and since $\mathbb{P}(D_i)\geq 1-\delta_i$ for $i=1,2$,
$$\mathbb{P}\left\{\xi(x)=\eta(x),\shs \forall\, x\in T_1\cup T_2\right\}\geq 1-(\delta_1+\delta_2).$$
We deduce by the same previous arguments that
$$\|\mathcal{P}-\mathcal{Q}\|_{var}\leq 4\,(\delta_1+\delta_2).$$

Finally, we have that
$$\|\mathcal{P}-\mathcal{P}_1\times\mathcal{P}_2\|_{var}\leq \|\mathcal{P}-\mathcal{Q}\|_{var}+\|\mathcal{Q}-\mathcal{Q}_1\times\mathcal{Q}_2\|_{var}+\|\mathcal{Q}_1\times\mathcal{Q}_2-\mathcal{P}_1\times\mathcal{P}_2\|_{var}.$$
As $\eta_1$ and $\eta_2$ are independent, 
$$\|\mathcal{Q}-\mathcal{Q}_1\times\mathcal{Q}_2\|_{var}=0.$$
Moreover,
$$\|\mathcal{P}_1\times\mathcal{P}_2-\mathcal{Q}_1\times\mathcal{Q}_2\|_{var}\leq\|\mathcal{P}_1-\mathcal{Q}_1\|_{var}+\|\mathcal{P}_2-\mathcal{Q}_2\|_{var}\leq 4\,(\delta_1+\delta_2)$$
and 
$$\|\mathcal{P}-\mathcal{Q}\|_{var}\leq 4\,(\delta_1+\delta_2).$$
Thus, we derive that
\begin{equation*}
\|\mathcal{P}-\mathcal{P}_1\times\mathcal{P}_2\|_{var}\leq 8\,(\delta_1+\delta_2).
\end{equation*}
\end{pro}

\subsubsection{Dimension $d=1$}

Remind that in this case $T_1=(-\infty,0]$ and $T_2=[r,+\infty)$ and denote by $\beta(r)$ the coefficient $\beta(T_1,T_2)$.

%------------------------------Theoreme 2-------------------------------------
\begin{thm}\label{reg1} (d=1) If Assumtions $1)$- $6)$ and $7a)$ are statisfied, the process $\xi$ has the absolute regularity property and for all $r>0$,
$$\beta(r)\leq C_1\textrm{e}^{-F(C_2 r)},$$
where the constants $C_1$ and $C_2$ can be choosen such that  $C_1=16$ and $C_2=\frac{1}{2M}.$
\end{thm}
We introduce for any subset $T$ of $\mathbb{R}$, the process $\xi_T$ defined as follows
\begin{equation}\label{xiT}
\xi_T(x)=\inf_{\scriptsize
\begin{array}{c}
g\in\mathcal{N}\\
x_g\in T
\end{array}}A_g(x) \hs\forall\, x\in\mathbb{R}^d.
\end{equation}
The proof of Theorem \ref{reg1} is based on the two following lemmas.
%-----------------------------Lemme 5-----------------------------------------
\begin{lem}\label{lemme1d1}
Under the same assumptions as in Theorem \ref{reg1}, for all $R>0$, we have that
$$\mathbb{P}\left\{\xi(x)=\xi_{(-\infty,MR]}(x),\shs\forall\,  x\leq 0\right\}\geq 1-\textrm{e}^{-F(R)}$$
with $\xi_{(-\infty,MR]}$ defined by relation (\ref{xiT}) with $T=(-\infty,MR]$.
\end{lem}

\begin{pro}Let us show first that
\begin{equation}\label{eq1}
\{\xi(0)\leq R\}\subset\{\xi(x)=\xi_{(-\infty,MR]}(x),\;\forall x\leq 0\}.
\end{equation}
Suppose that $\xi(0)\leq R$ and prove that
\begin{equation}\label{eq2}
\inf_{\scriptsize
\begin{array}{c}
g\in\mathcal{N}\\
x_g\leq MR
\end{array}}A_g(x)\leq \inf_{\scriptsize
\begin{array}{c}
g\in\mathcal{N}\\
x_g>MR
\end{array}}A_g(x)\hs\forall\, x\leq 0.
\end{equation}
For all $g=(x_g,t_g)\in E$ such that $x_g>MR$, assumptions $1)$ and $5)$ leads to
$$A_g(0)\geq t_g+\frac{|x_g|}{M}>R.$$
Since $\xi(0)\leq R$, we then deduce that
$$\xi(0)=\inf_{\scriptsize
\begin{array}{c}
g\in\mathcal{N}\\
x_g\leq MR
\end{array}}A_g(0).$$
Consequently, there exists $g_0\in\mathcal{N}$ such that $x_{g_0}\leq MR$ and $A_{g_0}(0)=\xi(0)$. Hence $A^-_g(0)>A_{g_0}(0)$ , for all $g=(x_g,t_g)\in \mathcal{N}$ such that $x_g>MR$ and  we deduce from Assumption $7a)$ that
$$A_g(x)\geq A_{g_0}(x)\hs\forall\, x\leq 0$$
and (\ref{eq2}) follows. Since 
$$\xi(x)=\min\{\inf_{\scriptsize
\begin{array}{c}
g\in\mathcal{N}\\
x_g\leq MR
\end{array}}A_g(x),\inf_{\scriptsize
\begin{array}{c}
g\in\mathcal{N}\\
x_g>MR
\end{array}}A_g(x)\},$$
we derive that 
$$\xi(x)=\xi_{(-\infty,MR]}(x)\hs\forall\, x\leq 0$$
and (\ref{eq1}) is then proved. Finally,
$$\mathbb{P}\left\{\xi(x)=\xi_{(-\infty,MR]}(x),\;\forall  x\leq 0\right\}\geq \mathbb{P}\left\{\xi(0)\leq R\right\}$$
and \begin{equation*}
\mathbb{P}\left\{\xi(0)\leq R\right\}\geq 1-\textrm{e}^{-\Lambda(K_{0,R})}.
\end{equation*}
\end{pro}
%-----------------------------Lemme 6-----------------------------------------
Thanks to symmetry arguments, we derive the following lemma.
\begin{lem}\label{lemme2d1}
Under the same assumptions as in Theorem \ref{reg1}, for all $R>0$, we have that
$$\mathbb{P}\left\{\xi(x)=\xi_{[MR,+\infty)}(x),\shs\forall\,  x\geq 2MR\right\}\geq 1-\textrm{e}^{-F(R)}$$
where $\xi_{[MR,+\infty)}$ is defined by relation (\ref{xiT}) with $T=[MR,+\infty)$.
\end{lem}
We turn back to the demonstration of Theorem \ref{reg1}.

\begin{pro} Let $r$ be a positive real and consider $R$ such that $2MR=r$ with $M$ the constant of Assumption $5)$. Lemma \ref{lemme1d1} and Lemma \ref{lemme2d1} allow us to apply Lemma \ref{outil} with $\eta_1=\xi_{(-\infty,MR]}$, $\eta_2=\xi_{[MR,+\infty)}$, $T_1=(-\infty,0]$, $T_2=[2MR,+\infty)$ and $\delta_1=\delta_2=e^{-F(R)}$. We obtain then that
\begin{equation*}
\beta(r)\leq 8\,(\delta_1+\delta_2)=16\,e^{-F(\frac{r}{2M})}.
\end{equation*}
\end{pro}

\subsubsection{Dimension $d\geq 2$}

We obtain first an upper bound for the absolute regularity coefficient in the case of two quadrants $T_1$ and $T_2$ which are separated by a $2r$-width band. As the random field $\xi$ is homogeneous, we can choose $T_1=\prod_{i=1}^d (-\infty,0]$ and $T_2=\prod_{i=1}^d[a_i,+\infty)$. We denote by $L_1$, (respectively $L_2$) the hyperplane orthogonal to $e=\frac{1}{\sqrt{d}}(1,\dots,1)$ and containing the point $(0,\dots,0)$ (respectively $(a_1,\dots,a_d)$) as represented on Figure $1$ when $d=2$. The distance between the hyperplanes $L_1$ and $L_2$ equals $2r=<e,a>$. Since $<e,a>$ is positive, we can introduce the hyperplane $L_0$ situated at equal distance between $L_1$ and $L_2$. Finally, we denote by $E_1$ (respectively $E_2$) the open half-space delimited by $L_0$ and containing $L_1$ (respectively $L_2$).

\begin{figure}[h]
\begin{center}
\label{figquad}
\includegraphics[width=6cm]{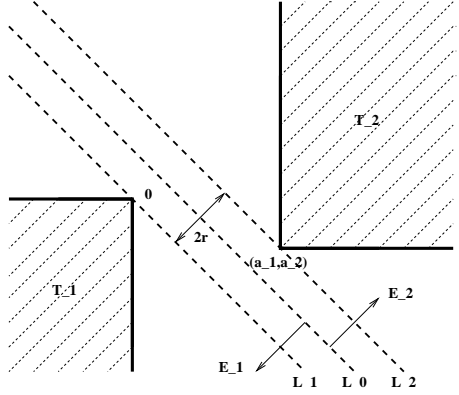}
\caption{Quadrant domains for $d=2$}
\end{center}
\end{figure}

%----------------------------------- Theoreme 3 -------------------------------
\begin{thm}\label{reg2} $(d\geq 2)$ If Assumptions $1)$- $9)$ are satisfied and $T_1$ and $T_2$ are the quadrant domains previously described, then
\begin{equation}\label{sigma1}
\beta(T_1,T_2)\leq 16\sum_{k=1}^{\infty}k^{d-1}\textrm{e}^{-F(C\,k)},
\end{equation}
where $F(t)$ is the measure of $K_t$, $C=\frac{2R}{dH}$, $R=\frac{r}{H}$ and $H=2(A+M)$ with $A$ and $M$ the constants of Assumptions $5)$ and $8)$.
\end{thm}
Before proving the theorem we give an estimate of the majorant series in  (\ref{sigma1}) for two typical cases. 

\begin{anex}\label{ex2}If $F(t)\geq (d+\delta)\ln t-\ln \gamma$ with $\delta,\,\gamma>0$, we have $\textrm{e}^{-F(t)}\leq \gamma \,t^{-(d+\delta)}$ and obtain a polynomial estimation of the sum:
$$\sum_{k=1}^{\infty}k^{d-1}\textrm{e}^{-F(Ck)}\leq \gamma'\,C^{-(d+\delta)}$$
with 
$$\gamma'=\gamma\,\sum_{k=1}^{\infty}k^{-(1+\delta)}.$$
\end{anex}
\begin{anex}\label{ex3}
If we rather suppose that $F(t)\geq \gamma \,t^{\delta}-c$ with $\delta,\,\gamma,\,c>0$, then $\textrm{e}^{-F(t)}\leq c_1\,\textrm{e}^{-\gamma \,t^{\delta}}$ with $c_1=\textrm{e}^c$. We derive a super-exponential estimation of the sum:
$$\sum_{k=1}^{\infty}k^{d-1}\textrm{e}^{-F(Ck)}\leq c_2\,\textrm{e}^{-\gamma \,C^{\delta}},$$
with $$c_2=c_1\sum_{k=1}^{\infty}k^{d-1}\textrm{e}^{-\gamma \,C^{\delta}\,(k^{\delta}-1)}.$$
\end{anex}
\begin{pro} 
Let us now introduce for all $r>0$ the following random fields:
\begin{equation}\label{eta1}
\eta_r^1(x)=\inf_{\begin{array}{c}
g\in\mathcal{N}\\
x_g\in E_1
\end{array}}A_g(x)\hs\forall\,x\in\mathbb{R}^d,
\end{equation}
\begin{equation}\label{eta2}
\eta_r^2(x)=\inf_{\begin{array}{c}
g\in\mathcal{N}\\
x_g\in E_2
\end{array}}A_g(x)\hs\forall\,x\in\mathbb{R}^d.
\end{equation}
For all $R>0$, we denote by $\xi_R$ the random field defined as follows:
\begin{equation}\label{xir}
\xi_R(x)=\inf_{\begin{array}{c}
g\in\mathcal{N}\\
|x_g|\leq R
\end{array}}A_g(x)\hs\forall\,x\in\mathbb{R}^d.
\end{equation}
The proof  is based on three lemmas.

%-------------------------------Lemme 7 ---------------------------------------
\begin{lem}\label{lemme1d2} Under the assumptions of Theorem \ref{reg2}, for all $R>0$,
\begin{equation*}
\mathbb{P}(\xi(x)=\xi_{HR}(x),\shs|x|\leq R)\geq 1-\textrm{e}^{-F(R)}
\end{equation*}
with $\xi_{HR}$ defined by equation (\ref{xir}).
\end{lem}
\begin{pro} We demonstrate that
$$\{\xi(0)\leq R\}\subset\{\xi(x)=\xi_{2(A+M)R}(x),\;|x|\leq R\}.$$
Suppose that $\xi(0)\leq R$ and consider $g_0\in\mathcal{N}$ such that $\xi(0)=A_{g_0}(0)$. Thanks to Assumption $5)$, we obtain that $|x_{g_0}|\leq MR$. On the other hand, we introduce $g=(0,R)$ and deduce from the definition of $\xi$ that
$$\xi(x)\leq A_{g_0}(x)\hs\forall\,x\in\mathbb{R}^d.$$
Since $g\in L_{g_0}$, Assumption $6)$ leads to
$$ A_{g_0}(x)\leq  A_{g}(x)\hs\forall\,x\in\mathbb{R}^d.$$
Thus,
$$\sup_{|x|\leq R}\xi(x)\leq \sup_{|x|\leq R}A_g(x).$$
Now for time 
$$t=\inf\{s>R\,|\,d_g(s)=2R\},$$ 
it is clear that
$$\sup_{|x|\leq R}A_g(x)\leq t$$
and from Assumption $8)$ we deduce that $D_g(t)\leq 2AR.$

Let us consider now a germ $g_1$ such that $|x_{g_1}|>2(A+M)R$. If $t_{g_1}\geq t$, Assumption $2)$ implies that 
$$A_{g_1}(x)\geq t\hs\forall\,x\in\mathbb{R}^d,\shs|x|\leq R.$$
If $t_{g_1}<t$, then Assumptions $3)$, $8)$ and $9)$ lead to
$$D_{g_1}(t)\leq D_{g_1}(R)+D_g(t).$$
But, assumption 5) implies that $D_{g_1}(R)\leq 2MR$ and the time $t$ is such that $D_g(t)\leq 2AR.$ Consequently,
$$D_{g_1}(t)\leq 2(A+M)R$$
and the crystals $C_g(t)$ and $C_{g_1}(t)$ are disjoint. Thus,
$$A_{g_1}(x)\geq t\hs\forall\,x\in\mathbb{R}^d,\shs|x|\leq R.$$
So, we conclude that 
\begin{equation*}
\xi(x)=\xi_{2(A+M)R}(x)\hs\forall\,x\in\mathbb{R}^d,\shs|x|\leq R.
\end{equation*}
\end{pro}

%------------------------------Lemme 8-----------------------------------------
\begin{lem} \label{lemme2d2}Under the assumptions of Theorem \ref{reg2}, for all $r>0$,
$$\mathbb{P}\{\xi(x)=\eta_r^1(x),\;x\in T_1\}\geq 1-\sum_{m=1}^{\infty}m^{d-1}\textrm{e}^{-F(Cm)}$$
with $\eta_r^1$ defined by (\ref{eta1}), $C=\frac{2R}{dH}$, $R=\frac{r}{H}$, $H=2(A+M)$.
\end{lem}
\begin{pro} We split the set $T_1$ into d-dimentional cubes denoted by $D_{\overline k}$, where for all $\overline k=(k_1,\dots,k_d)\in(-\mathbb{N})^d$,
$$D_{\overline k}=\prod_{i=1}^d[\frac{2R}{\sqrt d}(k_i-1),\frac{2R}{\sqrt d}(k_i)].$$
Each cube $D_{\overline k}$ is centered in $x_{\overline k}=(\frac{R}{\sqrt d}\left(2k_i-1)\right)_{i=1\dots d}$ and has diameter equal to $2R$. Remark also that the distance between $x_{\overline k}$ and $L_1$ equals $l_{\overline k}$ with
$$l_{\overline k}=R+|<\frac{2R}{\sqrt d}\overline{k},e>|=R(1+\frac{2}{d}|\sum_{i=1}^dk_i|).$$
Denote by $p$ the probability $\mathbb{P}(\xi(x)=\eta_r(x),\;x\in T_1)$ and note that
\begin{equation}\label{interbk}
p=\mathbb{P}(\bigcap_{\overline{k}\in(-\mathbb{N})^d} B_{\overline{k}}),
\end{equation}
with
$$B_{\overline{k}}=\{(\xi(x)=\eta_r(x),\;x\in D_{\overline k}\}.$$
From Lemma \ref{lemme1d2}, we obtain for all $a>0$ that
$$\mathbb{P}(\xi(x)=\xi_{B(x_{\overline k},Ha)}(x),\;|x-x_{\overline k}|\leq a)\geq 1-\textrm{e}^{-F(a)},$$
where
$$\xi_{B(x_{\overline k},Ha)}(x)=\inf_{\scriptsize\begin{array}{c}
g\in\mathcal{N}\\
|x_g-x_{\overline k}|\leq aH
\end{array}}A_g(x).$$
We choose $a=R+\frac{l_{\overline k}}{H}$. Hence, $D_{\overline k}\subset B(x_{\overline k},a)$ and
$$\{\xi(x)=\xi_{B(x_{\overline k},Ha)}(x),\;|x-x_{\overline k}|\leq a\}\subset\{\xi(x)=\xi_{B(x_{\overline k},Ha)}(x),\;x\in D_{\overline k}\}.$$
Moreover $B(x_{\overline k},Ha)$ is included in the half-space $E_1$. Consequently,
$$\{\xi(x)=\xi_{B(x_{\overline k},Ha)}(x),\;x\in D_{\overline k}\}\subset\{\xi(x)=\eta_r^1(x),\;x\in  D_{\overline k}\}.$$
Denoting by $p_{\overline k}$ the probability $\mathbb{P}(B_{\overline k})$, we finally obtain that
\begin{equation}\label{pk}
p_k\geq 1-\textrm{e}^{-F\left(R+\frac{l_{\overline k}}{H}\right)}.
\end{equation}
On the other hand, equation  (\ref{interbk}) implies that
\begin{equation*}
p=1-\mathbb{P}(\bigcup_{\overline k\in(-\mathbb{N})^d}B_{\overline{k}}^c).
\end{equation*}
From (\ref{pk}), we deduce that
\begin{equation}\label{somme}
p\geq 1-\sum_{\overline k\in(-\mathbb{N})^d}\textrm{e}^{-F\left(R+\frac{l_{\overline k}}{H}\right)}.
\end{equation}
Now, we obtain an upper bound for the sum in (\ref{somme}) as follows:
$$\begin{array}{ll}
 & \sum_{\overline k\in(-\mathbb{N})^d}\textrm{e}^{-F\left(R+\frac{l_{\overline k}}{H}\right)}\\
 & \\
 = & \sum_{m=0}^{\infty}\#\{\overline k,\;|\sum_{i=1}^d k_i|=m\}\,\textrm{e}^{-F\left(R+\frac{R}{H}(1+\frac{2}{d}m)\right)}\\
 & \\
 \leq & \sum_{m=0}^{\infty}(m+1)^{d-1}\textrm{e}^{-F\left(R\left(1+\frac{1}{H}(1+\frac{2}{d}m)\right)\right)}.\\
\end{array}$$
Since $R\left(1+\frac{1}{H}(1+\frac{2}{d}m)\right)\geq C(m+1)$ with $C=\frac{2R}{dH}$ when $d\geq 2$, we finally derive that
\begin{equation*}
p\geq 1-\sum_{m=1}^{\infty}m^{d-1}\textrm{e}^{-F(Cm)}.
\end{equation*}
\end{pro}
Symmetry arguments lead to the following lemma.
%---------------------------------Lemme 9--------------------------------------
\begin{lem}\label{lemme3d2}Under the assumptions of Theorem \ref{reg2}, for all $r>0$,
$$\mathbb{P}(\xi(x)=\eta_r^2(x),\;x\in T_2)\geq 1-\sum_{m=1}^{\infty}m^{d-1}\textrm{e}^{-F(Cm)}$$
with $\eta_r^2$ defined by (\ref{eta2}), $C=\frac{2R}{dH}$, $R=\frac{r}{H}$, $H=2(A+M)$.
\end{lem}
We make use of these three lemma to finish the proof of  Theorem \ref{reg2}.
 We apply Lemma \ref{outil} thanks to Lemma \ref{lemme2d2} and Lemma \ref{lemme3d2} with $\eta_1=\eta_1^r$, $\eta_2=\eta_2^r$, $\delta_1=\delta_2=\sum_{m=1}^{\infty}m^{d-1}e^{-F(C\,m)}$ and $T_1$, $T_2$ the quadrant domains. We then have that
\begin{equation*}
\beta(T_1,T_2)\leq 8\,(\delta_1+\delta_2)=16 \sum_{m=1}^{\infty}m^{d-1}e^{-F(C\,m)}.
\end{equation*}
\end{pro}

We give now an upper bound for the absolute regularity coefficient $\beta(T_1,T_2)$ in the case of enclosed cube domains separated by a $2r$-width polygonal band. As the random field $\xi$ is homogeneous, we consider centered domains $T_1=[-a,a]^d$ and $T_2=([-b,b]^d)^c$ as represented on Figure $2$ for $d=2$.

\begin{figure}[h]
\begin{center}
\label{figsqua}
\includegraphics[width=6cm]{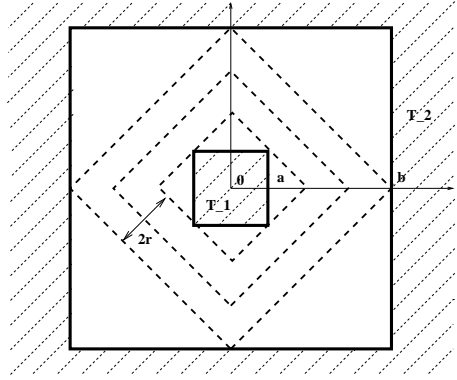}
\caption{Sketch for $d=2$}
\end{center}
\end{figure}
%-------------------------------Th4------------------------------
\begin{thm}\label{reg3} $(d\geq 2)$ If Assumptions $1)$-$9)$ are satisfied and $T_1$, $T_2$ are the enclosed domains previously described with $b\geq 2\,(H-1)\,a$, then
$$\beta(T_1,T_2)\leq 8\,(1+d\,2^d)\sum_{k=1}^{\infty}k^{d-1}\textrm{e}^{-F(C\,k)},$$
where $F(t)$ is the measure of $K_t$, $C=\frac{2R}{dH}$, $R=\frac{r}{H}$ and $H=2(A+M)$ with $A$ and $M$ the constants of Assumptions $5)$ and $8)$.
\end{thm}
The proof of Theorem \ref{reg3} make use of the same kind of arguments as in the proof of Theorem \ref{reg2}. Therefore, we intoduce first sets in order to define the random fields $\eta_1^r$ and $\eta_2^r$ approximating $\xi$ respectively on $T_1$ and $T_2$. Thus, we denote by $e_1,\dots,e_d$ the $d$ vectors of the canonical base in $\mathbb{R}^d$ and consider the set $A=\{\alpha\,|\,(\alpha_1,\dots,\alpha_d),\shs \alpha_i=\pm 1\}$ which cardinal equals $\# A=2^d$. For all $i$, the hyperplane $e_i^{\perp}$ separates the set $\mathbb{R}^d$ into two open half-space $E_i^{\epsilon}$ with $\epsilon=\pm 1$ and $\epsilon e_i$ contained in $E_i^{\epsilon}$. For all $\alpha\in A$, we introduce the quadrant:
\begin{equation*}
\mathcal{Z}_{\alpha}=\displaystyle\bigcap_{i=1}^d E_i^{\alpha_i}
\end{equation*} 
and for all $i=1\dots d$ the translated quadrant:
\begin{equation}\label{Zalpha}
\mathcal{Z}_{\alpha,i}=\mathcal{Z}_{\alpha}\varoplus \alpha_i\,b\,e_i.
\end{equation}
Observe that
\begin{equation*}
T_2=\displaystyle\bigcup_{\alpha\in A}\bigcup_{i=1}^d\mathcal{Z}_{\alpha,i}.
\end{equation*}
On the other hand, let us define for all $\alpha\in A$, the normed vector of $\mathcal{Z}_{\alpha}$:
$$d_{\alpha}=\frac{1}{\sqrt{d}}\sum_{i=1}^d\alpha_i\,e_i.$$
To separate the sets $T_1$ and $T_2$ by a $2r$-width polygonal band, the quantity $r=\frac{(b-2a)\sqrt{d}}{4}$ must be positive. Thus, we assume that $b\geq 2a$. In this case, we consider the hyperplanes 
$$\begin{array}{lcl}
L_{\alpha}^0 & = & d_{\alpha}^{\perp}+\frac{(b+2a)\sqrt{d}}{4}\,d_{\alpha}\\
 & & \\
L_{\alpha}^2 & = & L_{\alpha}^0 +r\,d_{\alpha}= d_{\alpha}^{\perp}+ \frac{b}{2}\sqrt{d}\,d_{\alpha}\\
 & & \\
L_{\alpha}^1  & = & L_{\alpha}^0 -r\,d_{\alpha}= d_{\alpha}^{\perp}+ a\,\sqrt{d}\,d_{\alpha}\\
\end{array}$$
as represented on Figure $3$ for $d=2$ and $\alpha=(1,1)$.

\begin{figure}[h]
\begin{center}
\label{figsqua}
\includegraphics[width=6cm]{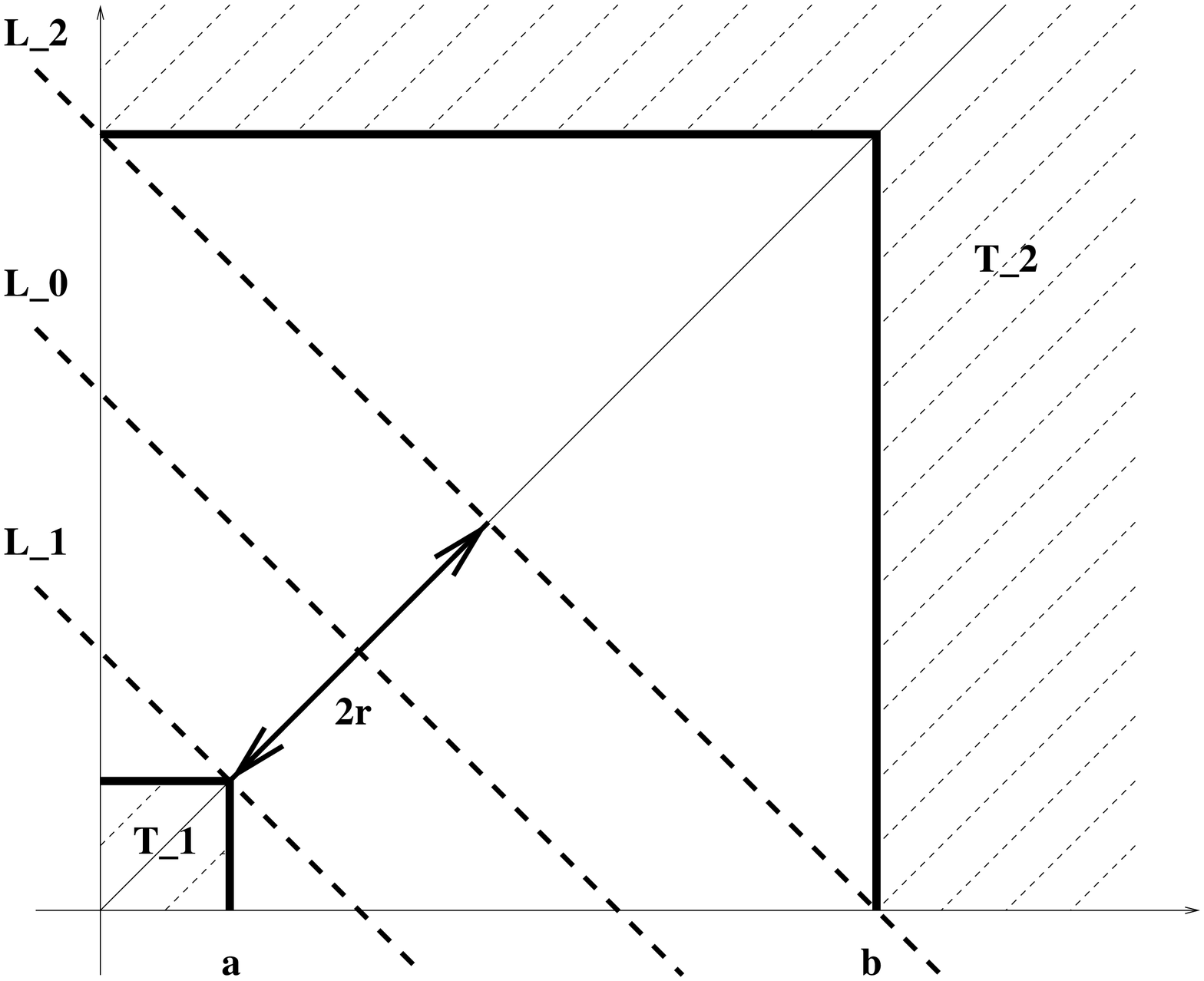}
\caption{Sketch for $d=2$}
\end{center}
\end{figure}

\noindent We introduce now for all $\alpha$ in $A$ the open half-space $S_{\alpha}^2$ delimited by the hyperplane $L_{\alpha}^0$ and containing the quadrants $\mathcal{Z}_{\alpha,i}$ for $i=1\dots d$. At last, we consider the set $S_2$ containing $T_2$:
\begin{equation*}
S_2=\displaystyle\bigcap_{\alpha\in A}S_{\alpha}^2.
\end{equation*}
Then, we introduce for all $\alpha\in A$, the random field:
$$\eta_{\alpha}(x)=\inf_{g\in S_{\alpha}^2}A_g(x)\hs\forall \, x\in\mathbb{R}^d$$
and approximate $\xi$ on $T_2$ by the following random field:
\begin{equation}\label{eta2bis}
\eta_2^r(x)=\inf_{g\in S_2}A_g(x)\hs\forall \, x\in\mathbb{R}^d.
\end{equation}

%-----------------------------------Lemme 10-----------------------------------
\begin{lem}\label{lemme4d2}If Assumptions $1)$-$9)$ are satisfied, then
$$\mathbb{P}\left\{\xi(x)=\eta_2^r(x),\shs x\in T_2\right\}\geq 1-d\,2^d\sum_{m=1}^{\infty}m^{d-1}e^{-F(C\,m)}$$
where $C$ is the constant of Theorem \ref{reg3} and $\eta_2^r$ is defined by (\ref{eta2bis}).
\end{lem}
\begin{pro} As for all $\alpha\in A$ and all $i=1\dots d$ the sets $\mathcal{Z}_{\alpha,i}$ defined by relation (\ref{Zalpha}) are quadrants included in $S_{\alpha}^2$, $\xi$ can be approximate by $\eta_{\alpha}$ on each $\mathcal{Z}_{\alpha,i}$ by Lemma \ref{lemme2d2} so that:
$$\mathbb{P}\left\{\xi(x)=\eta_{\alpha}(x),\shs\forall\,x\in \mathcal{Z}_{\alpha,i}\right\}\geq 1-\sum_{m=1}^{\infty}m^{d-1}e^{-F(C\,m)}.$$
Since,
$$\xi(x)\leq\eta_2^r(x)\leq \eta_{\alpha}(x),\hs\forall\,x\in\mathbb{R}^d$$
we deduce for all $\alpha\in A$ and all $i=1\dots d$ that
$$\mathbb{P}\left\{\xi(x)=\eta_2^r(x),\shs\forall\,x\in \mathcal{Z}_{\alpha,i}\right\}\geq 1-\sum_{m=1}^{\infty}m^{d-1}e^{-F(C\,m)}.$$
Finally, we derive that
\begin{equation*}
\mathbb{P}\left\{\xi(x)=\eta_2^r(x),\shs x\in T_2\right\}\geq 1-d\,2^d\sum_{m=1}^{\infty}m^{d-1}e^{-F(C\,m)}.
\end{equation*}
\end{pro}

We consider now for all $\alpha$ in $A$, the open half-space $S_{\alpha}^1=(S_{\alpha}^2)^c\backslash L_{\alpha}^0$. We also introduce the intersection
$$S_1=\displaystyle\bigcap_{\alpha\in A}S_{\alpha}^1$$
on which $\xi$ can be approximated by the following random field:
\begin{equation}\label{eta1bis}
\eta_1^r(x)=\inf_{g\in S_1}A_g(x)\hs\forall\,x\in\mathbb{R}^d.
\end{equation}
%-------------------------------Lemme 11--------------------------------------
\begin{lem}\label{lemme5d2}If Assumptions $1)$-$9)$ are satisfied with sets $T_1$ and $T_2$ such that $b\geq 2\,(H-1)\,a$, then
$$\mathbb{P}\left\{\xi(x)=\eta_1^r(x),\shs x\in T_1\right\}\geq 1-\sum_{m=1}^{\infty}m^{d-1}e^{-F(C\,m)}$$
where $C$ is the constant of Theorem \ref{reg3} and $\eta_1^r$ is defined by (\ref{eta1bis}).
\end{lem}

\begin{pro} We consider the centered open ball $B_1=B(0, a\,\sqrt{d})$ included in $T_1$ and the ball $B_2=B(0,a\,\sqrt{d}+r')$ with $r'\leq r$ so that $B_2$ is contained in $S_1$. If we denote by $R$ the radius of $B_1$ and assume that $R\,H=a\,\sqrt{d}+r'$ with $H$ the constant of Theorem \ref{reg3}, we derive that 
$$r'=(H-1)\,a\,\sqrt{d}\leq r=\frac{(b-2\,a)\sqrt{d}}{4}$$
 and finally that $b$ must be such that $b\geq 2\,(H-1)\,a$. Since $A\geq 1$, it follows that $H\geq 2$ and $b\geq2\,a$. We intoduce the random fields $\eta_{B_2}$:
\begin{equation*}
\eta_{B_2}(x)=\inf_{g\in B_2}A_g(x)\hs\forall\,x\in\mathbb{R}^d.
\end{equation*}
We remark that
\begin{equation}\label{majo}
\xi(x)\leq \eta_1^r(x)\leq \eta_{B_2}(x)\hs\forall\,x\in\mathbb{R}^d.
\end{equation}
But by Lemma \ref{lemme1d2},
$$\mathbb{P}\left\{\xi(x)=\eta_{B_2}(x),\shs\forall x\in B_1\right\}\geq 1-e^{-F(R)}$$
and from inequality (\ref{majo})
$$\mathbb{P}\left\{\xi(x)=\eta_1^r(x),\shs\forall x\in B_1\right\}\geq 1-e^{-F(R)}.$$
As $B_1\subset T_1$, we also have that
$$\left\{\xi(x)=\eta_1^r(x),\shs x\in B_1\right\}\subset\left\{\xi(x)=\eta_1^r(x),\shs x\in T_1\right\}$$
and then
$$\mathbb{P}\left\{\xi(x)=\eta_1^r(x),\shs\forall x\in T_1\right\}\geq 1-e^{-F(R)}.$$
Finally, as $H\geq 2$, $R\geq C$ with $C=\frac{2R}{dH}$ and $e^{-F(R)}\leq e^{-F(C)}$. We note that $e^{-F(C)}\leq \sum_{m=1}^{\infty}m^{d-1}e^{-F(Cm)}$, hence
\begin{equation*}
\mathbb{P}\left\{\xi(x)=\eta_1^r(x),\shs\forall x\in T_1\right\}\geq 1-\sum_{m=1}^{\infty}m^{d-1}e^{-F(Cm)}.
\end{equation*}
\end{pro}

\begin{pro}[Theorem \ref{reg3}] We apply again Lemma \ref{outil} thanks to Lemma \ref{lemme4d2} and Lemma \ref{lemme5d2} with $\eta_1=\eta_1^r$, $\eta_2=\eta_2^r$, $\delta_1=\sum_{m=1}^{\infty}m^{d-1}e^{-F(C\,m)}$, $\delta_2=d\,2^d\,\delta_1$ and $T_1$, $T_2$ the enclosed domains. We then have that
\begin{equation*}
\beta(T_1,T_2)\leq 8\,(\delta_1+\delta_2)=8\,(1+d\,2^d) \sum_{m=1}^{\infty}m^{d-1}e^{-F(C\,m)}.
\end{equation*}
\end{pro}

\subsection{Lower bounds}
In conclusion we give a lower bound of $\beta$-coefficient in the context of 
Examples \ref{ex2} and \ref{ex3} which are of the same type as the upper ones.
It shows that
the upper bounds in Theorem \ref{reg1}, Theorem \ref{reg2} and Theorem \ref{reg3} are sufficiently precise.

Let the dimension $d =1$. We choose $A=\{\xi(0)>a\}$ and $B=\{\xi(x)>a\}$ with $|x|=r$. It is clear that
\begin{equation}\label{minor}
\beta(r)\geq 2|\mathbb{P}(A\cap B)-\mathbb{P}(A)\mathbb{P}(B)|.
\end{equation}

Since $\xi$ is space homogeneous, we obtain that
$$\begin{array}{lcl}
\mathbb{P}(A)=\mathbb{P}(B) & = & \mathbb{P}\left\{\mathcal{N}\cap K_a=\emptyset\right\}\\
 & = & e^{-F(a)}.
\end{array}$$
To compute $\mathbb{P}(A\cap B)$, we assume that there exists $\tau>0$ such that for all $t$ large enough,
$$K_t\cap(K_t+h)\subset K_{(1+\tau)\,t}\hs\forall\,h\in\mathbb{R}^d,\shs|h|\leq \tau t.$$
Under this assumption, 
$$\begin{array}{lcl}
\mathbb{P}(A\cap B) & = & \mathbb{P}\left\{\mathcal{N}\cap K_a=\emptyset,\shs \mathcal{N}\cap(K_a+x)=\emptyset\right\}\\
 & \geq & \mathbb{P}\left\{\mathcal{N}\cap K_{(1+\tau)\,a}=\emptyset\right\}\\
 & = & e^{-F\left((1+\tau)\,a\right)}.
\end{array}$$
We choose, $r=\delta\,a$ so that 
\begin{equation}\label{betaminor}
\beta(r)\geq\left|e^{-2\,F\left(\frac{r}{\tau}\right)}-e^{-F\left(\frac{(1+\tau)}{\tau}r\right)}\right|
\end{equation}
We compute the minoration term in inequality (\ref{betaminor}) for the two examples. In the case of Example \ref{ex2} where $F(t)=(d+\delta)\,\ln(t)-\ln(\gamma)$ with $\delta,\,\gamma>0$, we obtain that
$$e^{-2\,F\left(\frac{r}{\tau}\right)}=\gamma^2\tau^{2\,(d+\delta)}r^{-2\,(d+\delta)}$$
and 
$$e^{-F\left(\frac{(1+\tau)}{\tau}r\right)}=\gamma\left(\frac{\tau}{\tau+1}\right)^{d+\delta}r^{-(d+\delta)}.$$
Thus, for $r$ sufficiently large,
$$\beta(r)\geq \kappa_1r^{-(d+\delta)}$$
with $\kappa_1>0$. For Example \ref{ex3} where $F(t)=\gamma\,t^\delta-c$ with $\gamma,\,\delta,_,c>0$, we derive that
$$e^{-2\,F\left(\frac{r}{\tau}\right)}=e^{2\,c}e^{-\frac{2\,\gamma}{\tau^\delta}r^\delta}$$
and
$$e^{-F\left(\frac{(1+\tau)}{\tau}r\right)}=e^ce^{-\frac{\gamma\,(1+\tau)^\delta}{\tau^\delta}r^\delta}.$$
Finally, if $\tau<2^{\frac{1}{\delta}}-1$, then for $r$ sufficiently large,
$$\beta(r)\geq \kappa_2 e^{-\gamma\left(\frac{1+\tau}{\tau}\right)^\delta r^\delta}$$
with $\kappa_2>0.$

%\bibliographystyle{plain}
%\bibliography{ref}

\begin{thebibliography}{00}
% please try to use the bibitem system -
% the references should be in alphabetical order of authors' names.
% Articles with a single author first, author will 1 co-author next,
% then author with several co-authors;




%\bibitem{DI06}  Yu. Davydov and A. Illig,  Ergodic properties of  geometrical cristallization processes,  CRAS (2007)


\bibitem{JM39}  W. A. Johnson, and R. F. Mehl, Reaction Kinetics in Processes of Nucleation and Growth, Trans. Amer. Inst. Min. Metal. Petro. Eng. 135 (1939), pp. 416-458.

\bibitem{MC97}  A. Micheletti, and V. Capasso, The stochastic geometry of polymer crystallization processes, Stochastic Anal. Appl. 15 (1997) no. 3, pp. 355-373.


\bibitem{Kol37} A. N. Kolmogorov, Statistical theory of crystallization of metals, Bull. Acad. Sci. USSR Mat. Ser. 1 (1937) pp. 355-359.

\bibitem{Mol86} J. M{\o}ller, Random tessellations in {${\bf R}\sp d$}, Memoirs of Aarhus University Institute of Mathematics Department of Theoretical Statistics 9 (1986).

\bibitem{Mol89} J. M{\o}ller, Random tessellations in {${\bf R}\sp d$},  Adv. in Appl. Probab. 21 (1989) pp. 37-73.


\bibitem{Mol92} J. M{\o}ller, Random {J}ohnson-{M}ehl tessellations, Adv. in Appl. Probab. 24 (1992) pp. 814-844. 

\bibitem{Mol95} J. M{\o}ller, Generation of {J}ohnson-{M}ehl crystals and 
comparative analysis of models for random nucleation, Adv. in Appl. Probab. 27 (1995) pp. 367-383.







\end{thebibliography}

\end{document}